# SPECIAL INVITED PAPER

## COEXISTENCE IN STOCHASTIC SPATIAL MODELS[1]

By Rick Durrett

*Cornell University*

In this paper I will review twenty years of work on the question: When is there coexistence in stochastic spatial models? The answer, announced in Durrett and Levin [*Theor. Pop. Biol.* **46** (1994) 363–394], and that we explain in this paper is that this can be determined by examining the mean-field ODE. There are a number of rigorous results in support of this picture, but we will state nine challenging and important open problems, most of which date from the 1990's.

**Introduction.** There is an incredible diversity of species that coexist in the world. At the Botanic Garden in Singapore one can see 1000 species of orchids. These are cultivated, of course, but if one examines the food web in a small lake one finds dozens of species coexisting. An important problem in ecology is to identify mechanisms that permit the coexistence of species. In this paper we will examine that question in the context of stochastic spatial models. In these interacting particle systems, space is represented by the $d$-dimensional integer lattice $\mathbb{Z}^d$. With ecological problems in mind, we will usually take $d = 2$.

Our story begins with a very simple and natural model, but one that is still not well understood.

*Example* 0. *Competing contact processes.*

- Each site in $\mathbb{Z}^2$ can be in state $0 =$ vacant, or in state $i = 1, 2$ to indicate that it is occupied by one individual of type $i$.

Received August 2008; revised December 2008.

[1]Based on the third Wald lecture given at the World Congress of Probability and Statistics, Singapore, July 14–19, 2008. Supported in part by an NSF grant from the probability program.

*AMS 2000 subject classification.* 60K35.

*Key words and phrases.* Interacting particle system, block construction, fast stirring limit, competitive exclusion principle.

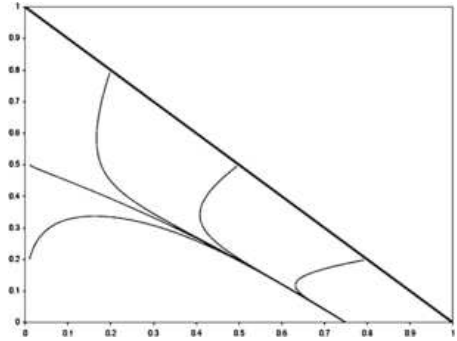

Fig. 1. *Competing contact process ODE.*

- Individuals of type $i$ die at rate $\delta_i$ and give birth at rate $\beta_i$. Here, at rate $\lambda$ means that these events happen at times of a rate $\lambda$ Poisson process.
- A type $i$ born at $x$ goes to $x + y$ with probability $p_i(y)$. If $x + y$ is vacant it changes to state $i$, otherwise nothing happens.

When there is only one type this reduces to the system introduced by Harris (1974). After several decades of work, this model is very well understood. See Liggett (1999) for a survey.

If we assume that the states of adjacent sites are independent then the fraction of sites $u_i$ in state $i = 1, 2$ satisfies

$$
\begin{aligned}
\frac{du_1}{dt} &= \beta_1 u_1 (1 - u_1 - u_2) - \delta_1 u_1, \\
\frac{du_2}{dt} &= \beta_2 u_2 (1 - u_1 - u_2) - \delta_2 u_2.
\end{aligned}
\tag{1}
$$

This is called the *mean-field ODE*, because if we consider the system on $N$ sites with a uniform dispersal distribution then in the limit as $N \to \infty$ the densities converge to this limit. In the spatial model adjacent sites are not independent. However, writing and analyzing the mean-field ODE is a good first step in guessing what the system will do.

In (1) $du_i/dt = 0$ when $(1 - u_1 - u_2) = \delta_i/\beta_i$. These lines are parallel, so they either do not intersect or coincide. Figure 1 shows the mean-field ODE when $\beta_1 = 4$, $\beta_2 = 2$, and $\delta_1 = \delta_2 = 1$. In this case all solutions starting from a point with $u_1 > 0$ converge to $(3/4, 0)$.

Neuhauser (1992) proved the following result:

THEOREM 1. *If the dispersal distributions $p_1 = p_2 = p$, death rates $\delta_1 = \delta_2 = \delta$, and birth rates $\beta_1 > \beta_2$ then species 1 out competes species 2. That is, if the initial condition is translation invariant and has $P(\xi_0(x) = 1) > 0$ then $P(\xi_t(x) = 2) \to 0$.*

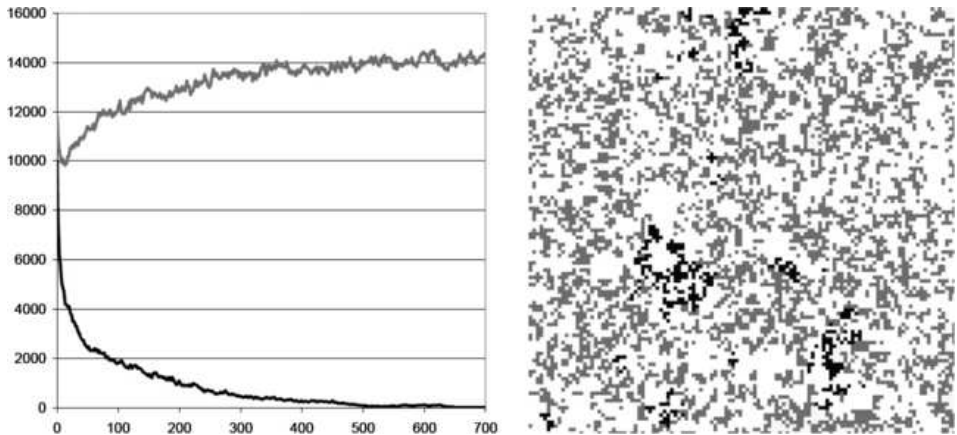

FIG. 2. *Simulation of competing contact process.* $\beta_1 = 3.9$, $\delta_1 = 2$ *(black) versus* $\beta_2 = 2$, $\delta_1 = 1$ *(gray). The picture is a snapshot of part of the grid at time 300.*

IDEAS BEHIND THE PROOF. We construct the process from a "graphical representation." For each site $x$ there is a rate $\delta$ Poisson process $D_n^x$, $n \geq 1$ that kills the particle at $x$ (if there is one). For each $x$ and $y$, there are Poisson processes $B_n^{x,y}$, $n \geq 1$ and $A_n^{x,y}$, $n \geq 1$ with rates $\beta_2 p(y)$ and $(\beta_1 - \beta_2) p(y)$. The first causes births from $x$ to $x + y$ if $x$ is occupied and $x + y$ is vacant. The second causes births $x$ to $x + y$ if $x$ is occupied by a 1 and $x + y$ is vacant. Using this construction and working backwards in time, Neuhauser (1992) was able to show that the extra arrows reserved only for the 1's meant that if the site was occupied then it would be a 1 with high probability. We refer the reader to Neuhauser (1992) for the details, which are somewhat complicated. Our reason for giving this sketch is to make clear that the proof requires $\delta_1 = \delta_2$. □

In the mean-field ODE, only the ratios $\beta_i / \delta_i$ matters, so it is natural to guess.

PROBLEM 1. Show that species 1 outcompetes species 2 holds if the dispersal distributions are the same and $\beta_1 / \delta_1 > \beta_2 / \delta_2$. The simulation in Figure 2 gives some support for this conjecture.

The behavior in Problem 1 is what biologists expect based on the *Competitive Exclusion Principle*. A version of this can be found in work of Levin (1970). Consider an ODE of the form:

$$\frac{du_i}{dt} = u_i f_i(z_1, \ldots, z_m), \qquad 1 \leq i \leq n.$$

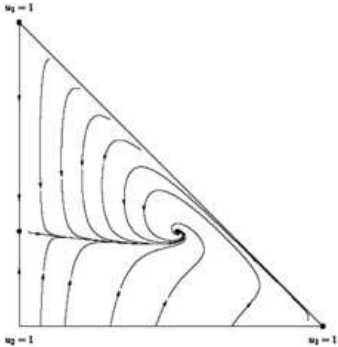

Fig. 3. *Host pathogen mean-field ODE.*

Here the $z_i$ are resources. In the previous model there is one resource: $z_1 = 1 - u_1 - u_2$ free space. A special case is shown in Figure 3. Omitting the details of the assumptions on the $f_i$ in the mathematical result, the principle as biologists use it is:

THEOREM 2. *If $n > m$ then no stable equilibrium in which all $n$ species are present is possible.*

In words, the number of coexisting species is smaller than the number of resources.

To try to build some excitement about Problem 1, we note that Chan and Durrett (2006) proved that a "fugitive species" that disperses at a fast rate, according to a truncated power law distribution, can coexist with a superior competitor with a nearest neighbor dispersal distribution. This does not contradict the competitive exclusion principle, because in addition to single site deaths their model has forest fires, which destroy large squares of occupied sites at a small rate. Thus the model has a second type of space, "recently disturbed space," and is entitled to have two coexisting species.

As announced in the abstract, one goal of this paper is to explain the idea of Durrett and Levin (1994) that one can determine whether coexistence happens in the stochastic spatial model by examining properties of the mean-field ODE. The discussion is divided into three cases according to the properties of the ODE.

**Case 1. Attracting fixed point.** When the mean-field ODE has an attracting fixed point with all components positive, we expect coexistence in the spatial model, that is, there is a stationary distribution which concentrates on configurations that have infinitely many sites occupied by each species.

*Example* 1.1. *Grass bushes trees.* In this variant of the contact process there is a hierarchy of types. In hindsight this is a very natural model. However, it owes its invention to talking to Simon Levin about successional sequences in a forest.

- Each site in $\mathbb{Z}^2$ can be in state $0 =$ grass, $1 =$ bush, $2 =$ tree.
- Individuals of type $i$ die at rate $\delta_i$ and give birth at rate $\beta_i$.
- A type $i$ born at $x$ goes to $x+y$ with probability $p_i(y)$. If $x+y$ is in state $j < i$ it changes to state $i$, otherwise nothing happens.

The mean field ODE is

$$
\begin{aligned}
\frac{du_1}{dt} &= \beta_1 u_1(1 - u_1 - u_2) - \delta_1 u_1 - \beta_2 u_2 u_1, \\
\frac{du_2}{dt} &= \beta_2 u_2(1 - u_1) - \delta_2 u_2.
\end{aligned}
\tag{2}
$$

If $\beta_2 > \delta_2$, the equilibrium frequency of 2's, $u_2^* = (\beta_2 - \delta_2)/\beta_2$. Given this, one can solve for $u_1^*$ and see when it is positive. However, it is better to approach the question by examining when the "1's can invade 2's in equilibrium." The phrase in quotes means that if the 2's are in equilibrium and the 1's are at a small density then the density of 1's will increase. Ignoring the possibility that a 1 will see another 1 nearby, the condition is:

$$
\beta_1 \cdot \frac{\delta_2}{\beta_2} > \delta_1 + \beta_2 \cdot \frac{\beta_2 - \delta_2}{\beta_2}.
\tag{3}
$$

The left-hand side gives the rate at which 1's give birth onto 0's, while on the right, the first term is the death rate of 1's and the second is the rate at which they are eliminated by births of 2's. It is easy to check that if $\beta_2 > \delta_2$ and (3) is satisfied then (2) has an equilibrium will all components positive.

When $\delta_1 = \delta_2 = 1$, (3) becomes $\beta_1 > \beta_2^2 > 1$. For simplicity, we will consider only this case. Durrett and Swindle (1991) have shown

THEOREM 3. *If $\beta_1 > \beta_2^2 > 1$ then when $p_i$ is uniform on $\{x : 0 < \|x\| \leq L\}$ and $L$ is large, there is a stationary distribution $\mu_{12}$ that concentrates on configurations with infinitely many 1's and 2's.*

OUR HAMMER: THE BLOCK CONSTRUCTION. The survival of the 2's is not a problem because they are a contact process and don't feel the presence of the 1's. To prove that the 1's can persist in the space that remains, we use a "block construction," which consists of comparing the particle system with a mildly dependent oriented percolation in which sites are open with probability close to 1.

For an account of this method see my St. Flour lecture notes, Durrett (1995), or return to the first application, Bramson and Durrett (1988), for a very simple example. At an intuitive level, what one shows is that "one pile will make two piles with high probability."

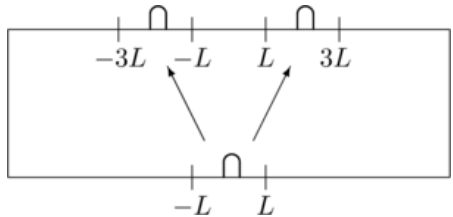

Here, pile is an undefined term and the phrase is short hand for the main idea behind the construction. To be more precise, in $d=1$ one might show that $L$ and $T$ can be chosen so that if $[-L, L]$ is "good" at time 0, then $[-3L, -L]$ and $[L, 3L]$ will be good at time $T$ with probability $1 - \varepsilon$. In this case, good might be there are not too many 2's in $[-2L, 2L]$ and there are enough 1's in $[-L, L]$. To get a finite range of dependence between the events in the construction, we need to estimate the probability of the good events if we assume that all sites in $[-kL, kL]^c$ are always occupied by 2's.

If $\varepsilon < \varepsilon_k$, that is, a constant that depends on the range of dependence, then facts about oriented percolation give a lower bound on the density of sites occupied by 1's. Taking the Cesaro average of the distribution from time 0 to time $t$, and finding a convergent subsequence produces the desired stationary distribution. For more details on the last point, see Liggett (1985). □

The block construction technology has improved quite a bit since 1991, so at this point it should be fairly routine to do the following:

EXERCISE. Show that if $\delta_1 = \delta_2 = 1$, $\beta_2 > 1$ and $\beta_1 < \beta_2^2$ then the 1's die out when the range is large.

In the setting of Theorem 3, in addition to existence of a stationary distribution, we have a uniqueness result proved by Durrett and Moller (1991).

THEOREM 4. *Suppose $\delta_1 = \delta_2 = 1$, $\beta_1 > \beta_2^2 > 1$. If the range is large then whenever the 1's and the 2's do not die out then the process converges to $\mu_{12}$.*

There are stationary distributions $\mu_1$ and $\mu_2$ with only 1's and 2's respectively. By results for the one-type contact process, these are unique if one specifies that there is no mass on the all 0's state. In addition there is the trivial stationary distribution $\mu_0$ that assigns mass 1 to all 0's. The convergence result in Theorem 4 when combined with results for the one-type contact process implies that all stationary distributions are convex combinations of $\mu_0$, $\mu_1$, $\mu_2$ and $\mu_{12}$.

*Interlude: fast stirring.* There are, at this point, a number of coexistence results for particle systems with long range interactions: Durrett (1992), Durrett and Schinazi (1993), Durrett and Neuhauser (1997), Durrett and Lanchier (2008), etc. However, the proofs of these results are done on a case by case basis. Things are simpler if, instead of long range, we assume that there is "fast stirring": suppose that for each pair of nearest neighbors $x$ and $y$, at rate $\varepsilon^{-2}$ exchange the values $\xi_t(x)$ and $\xi_t(y)$. In this case there is a general result.

THEOREM 5. *Suppose there is a function $\phi$ that* (i) *decreases along solutions of the mean-field ODE,* (ii) *is convex, and* (iii) $\phi(u) \to \infty$ *when* $\min_i u_i \to 0$. *Then there is coexistence in the model with fast stirring when* $\varepsilon < \varepsilon_0$.

Durrett (2002) applies this result to a wide variety of systems: epidemics, predator-prey models, predator mediated coexistence, etc.

SKETCH OF THE PROOF OF THEOREM 5. Suppose that the mean-field ODE is

$$\frac{du}{dt} = f(u).$$

Conditions (i)–(iii) and a few lines of calculus implies, see pages 102–103 of Durrett (2002) that if $h = \phi(u)$ and $\hat{\phi}$ is the time derivative of $\phi$ along solutions of the ODE then

$$\frac{dh}{dt} \leq \Delta h + \hat{\phi}(u).$$

Since $\hat{\phi}(u) \leq 0$, this is Brownian motion with killing. Using this we can show that there are constants $\delta, c > 0$ and $T < \infty$, so given starting conditions with $u_i(0, x) \geq \eta > 0$ for $x \in [-L, L]$, solutions of the "mean-field PDE"

$$\frac{du}{dt} = \Delta u + f(u)$$

have $\min_i u_i(t, x) \geq \delta$ for $t \geq T$, $|x| \leq ct$. As $\varepsilon \to 0$ the particle system on $\varepsilon \mathbb{Z}^d$ converges to the solution of the PDE. Using this with the result for the PDE, we have shown that "one pile will make two piles with high probability" and the result follows from the block construction. □

PROBLEM 2. Formulate and prove a similar general result for systems with long-range interactions.

*Example* 1.2. *Host-pathogen models.* In Durrett (2002) it was shown that predation can cause two competing species to coexist. Durrett and Lanchier (2008) have shown that coexistence can occur if there is a pathogen in one species. In the next models 1 and 3 are the two species, while 2 is species 1 in the presence of a pathogen. Letting $f_i$ be the fraction of neighbors in state $i$, the rates are

$$\begin{aligned} 1 &\to 2 & \alpha f_2 \\ 2 &\to 1 & \gamma_2(f_1 + f_2) \\ 3 &\to 1 & \gamma_3(f_1 + f_2) \\ 1 &\to 3 & \gamma_1 f_3 \\ 2 &\to 3 & \gamma_2 f_3 \end{aligned}$$

The first two rates say that if there are no 3's and we think of $1 = $ vacant and $2 = $ occupied then the 1's and 2's are a contact process. To explain the last four rates: at rate $\gamma_i$ individuals of type $i$ are replaced by the offspring of a randomly chosen neighbor. If the neighbor is type 1 or type 3 then the offspring has the same type as the parent. However, if the neighbor is type 2, the offspring is of type 1 because the pathogen is not passed into the seeds.

The mean-field ODE is

$$\begin{aligned} \frac{du_1}{dt} &= (u_1 + u_2)(\gamma_2 u_2 + \gamma_3 u_3) - \alpha u_1 u_2 - \gamma_1 u_1 u_3, \\ \frac{du_2}{dt} &= \alpha u_1 u_2 - \gamma_2 u_2, \\ \frac{du_3}{dt} &= u_3(\gamma_1 u_1 + \gamma_2 u_2) - \gamma_3 u_3(u_1 + u_2). \end{aligned} \tag{4}$$

As we will explain in a moment this leads easily to the following

THEOREM 6. *Suppose* $\gamma_1 < \gamma_3 < \gamma_2 < \alpha$ *and*

$$\gamma_1 \frac{\gamma_2}{\alpha} + \gamma_2 \left(1 - \frac{\gamma_2}{\alpha}\right) > \gamma_3 \tag{5}$$

*then there is coexistence for large range.*

KEYS TO THE PROOF. Here the 1's and 2's are a contact process so on the boundary $u_3 = 0$, $u_1 = \gamma_2/\alpha$ is an attracting fixed point. The displayed condition says that the 3's can invade the 1's and 2's in equilibrium. The key to the proof is using an understanding of the ODE to show that if the density of some type becomes small then a sequence of events will occur that results in all of the densities being bigger than some $\varepsilon$. The two observations which lead to this are: (i) the boundary $u_1 = 0$ is not a problem since 2's give birth to 1's; (ii) on the boundary $u_2 = 0$, 1's and 3's are a biased voter model in which 1's outcompete 3's. □

Figure 4 gives two simulations. It is clear that one cannot have coexistence if $\gamma_1$ and $\gamma_2$ are both $> \gamma_3$ or both $< \gamma_3$. The next problem address the remaining case:

PROBLEM 3. Show that coexistence is not possible in the host-pathogen model if $\gamma_2 < \gamma_3 < \gamma_1$.

In this case the pathogen is called a mutualist, since it decreases the rate at which the species is replaced. To explain why the result in Problem 3 should be true, note that if we start with the 1's and 2's in equilibrium and a small density of 3's then once the invasion of the 3's starts the fraction of 2's gets smaller, and the 3's have an even bigger advantage over the 1's and 2's. It is not hard to check in this case that there is no interior fixed point. For more results and problems about host-pathogen models see Lanchier and Neuhauser (2006).

**Case 2. Two locally attracting fixed points.** In this case, the limiting behavior of the ODE depends on the initial density. However, this is not the expected behavior for the particle system, and in this case the outcome of competition is dictated by the behavior of the PDE. The reason is that even if the initial distribution is translation invariant and hence has a well-defined density, there will be regions of space where the density of 1's is close to 1 and others where it is close to 0. To explain what we expect to happen in the particle system, we consider an example.

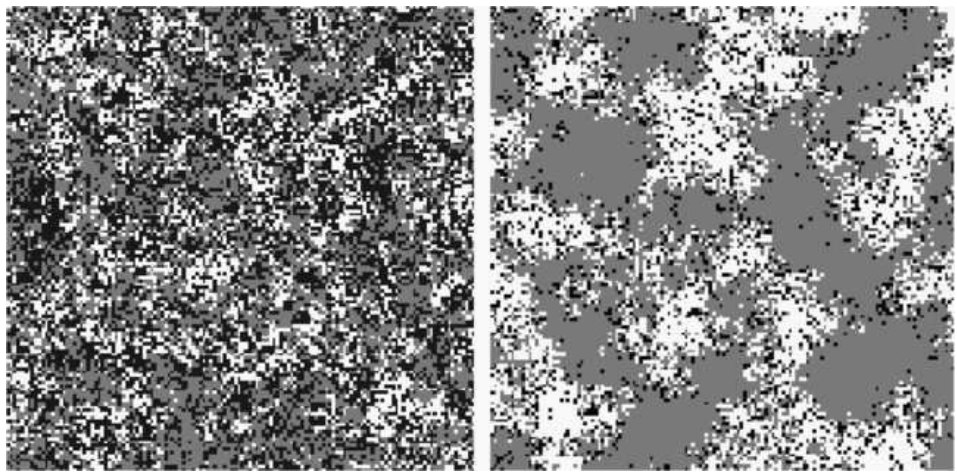

FIG. 4. *Host pathogen systems with coexistence and no coexistence. 1 = black, 2 = white, 3 = gray.*

*Example* 2.1. *Sexual reproduction.* The flip rates are as follows:

- $1 \to 0$ at rate 1,
- $0 \to 1$ at rate $\beta k(k-1)/n(n-1)$ if $k$ of the $n$ neighboring sites are occupied.

In words, at rate $\beta$, a vacant site picks two of its neighbors at random and become occupied if they both are.

The mean-field equation is:

$$\text{(6)} \qquad \frac{du}{dt} = -u + \beta u^2(1-u) = u(-1 + \beta u(1-u)).$$

Remembering that $u(1-u)$ is maximized at $1/2$, where the value is $1/4$, we see that there are nontrivial fixed points $\rho_1 < \rho_2$ if and only if $\beta > 4$, while if $\beta = 4$, $1/2$ is a double root.

At this point one might think that in the presence of fast stirring, the critical value for survival of the process $\beta_c \approx 4$ but that is not correct. To determine the asymptotics of the critical value you have to consider the mean-field PDE:

$$\text{(7)} \qquad \frac{\partial u}{\partial t} = \Delta u + g(u),$$

where $g(u) = u(-1 + \beta u(1-u))$.

A solution of (7) of the form $u(t,x) = w(x - ct)$ with $w(-\infty) = \rho_2$ and $w(+\infty) = 0$ is called a traveling wave. In order to be a solution $w$ must satisfy

$$-cw' = w'' + g(w).$$

Multiplying by $w'$ and integrating from $-\infty$ to $\infty$

$$-c \int (w')^2 \, dx = \int w'' w' \, dx + \int g(w) w' \, dx.$$

The first antiderivative on the right is $(w')^2/2$ which vanishes at $-\infty$ and $\infty$. Changing variables $y = w(x)$ in the second we have

$$c \int (w')^2 \, dx = \int_0^{\rho_2} g(y) \, dy.$$

Thus the sign of $c =$ the sign of $\int_0^{\rho_2} g(y) \, dy$. When $\beta = 4.5$, the cubic has roots at $2/3$, $1/3$, and $0$ and symmetry dictates that $\int_0^{\rho_2} g(y) \, dy = 0$. Thus $c > 0$ if and only if $\beta > 4.5$.

Durrett and Neuhauser (1994) have shown:

THEOREM 7. *If we introduce fast stirring at rate $\varepsilon^{-2}$ then in the sexual reproduction model, $\beta_c \to 4.5$ as $\varepsilon \to 0$.*

SKETCH OF THE PROOF. The key is the PDE fact that if $\eta > 0$ and the initial condition is $u(0,x) \geq \rho_1 + \eta$ for $|x| \leq L$ and $L \geq L_\eta$ then for any $\delta > 0$ we have $u(t,x) \approx \rho_2$ for $|x| \leq (c-\delta)t$ where $c$ is the wave speed defined above. Combining this with the convergence of the particle system with fast stirring to the PDE, we have the "one pile makes two piles with high probability" needed for the block construction. □

There are many situations in which we can prove the existence of stationary distributions, but convergence results like the one for grass-bushes-trees in Theorem 4 are rare. One exception is the result for multicolor contact processes in Durrett (1992). One would hope to develop general methods for proving uniqueness, but it seems sensible to start with a concrete case.

PROBLEM 4. Consider the sexual reproduction model with $\beta > 4.5$. Show that when the stirring rate $\varepsilon^{-2}$ is large there is a unique nontrivial stationary distribution.

*Example 2.2. Catalyst.* Moving away from ecology, our next system is a model for the catalytic converter in a car's exhaust system. States are $0 =$ vacant, $1 =$ CO (carbon monoxide attached to the surface), $2 =$ O (single oxygen atom attached to the surface). The rates are as follows:

- $0 \to 1$ at rate $p$,
- a pair of neighboring 0's $\to 22$ at rate $q/4$,
- adjacent $12 \to 00$ at rate $r/4$ (reaction to form $CO_2$).

In this model all 1's and all 2's are absorbing states corresponding to poisoning of the catalyst surface. In order for the catalytic converter to work and turn CO into $CO_2$ there must be coexistence in the spatial model. Ziff, Gulari and Barshad (1986) considered the case in which $r = \infty$ and $q/2 = 1 - p$ (the latter condition can be imposed by scaling time). Their simulations shows coexistence for $0.389 \leq p \leq 0.525$.

Since 1's land at rate $p$ and two 2's land at rate $\leq q = 2(1-p)$, it seems clear that the system converges to all 1's when $p \geq q$. There is a simple argument, see Theorem 1 in Durrett and Swindle (1994), which shows that if $p \geq q$ then $P(\xi_t(x) = 0) \to 0$ and if $x$ and $y$ are neighbors, $P(\xi_t(x) = 1, \xi_t(y) = 2) \to 0$, so

$$P(\xi_t(x) \equiv 1 \text{ on } [-K,K]^2) + P(\xi_t(x) \equiv 2 \text{ on } [-K,K]^2) \to 1$$

but we do not know how to prove that if we start from the $\equiv 0$ configuration, the system converges to all 1's. It is not hard to show, see Theorem 2 in Durrett and Swindle (1994), that the system converges to all 2's for small $p$. However, it is much more interesting to PROBLEM 5. Find $0 < p_1 < p_2 < 1$ so that coexistence occurs in the catalyst model for $p \in (p_1, p_2)$.

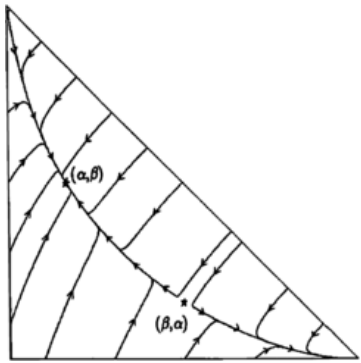

Fig. 5.  *Mean-field ODE for the catalyst.*

Simulations suggest that the density of O atoms in equilibrium drops to 0 discontinuously at the upper critical value, but proving this is a much harder problem.

Bramson and Neuhauser (1992) have proved coexistence when the $O_2$'s are replaced by an $N \times N$ polymer, which needs a vacant $N \times N$ square to land, and $N$ is sufficiently large. Durrett and Swindle (1994) proved coexistence in the original model by introducing fast stirring. The mean-field PDE is:

(8)
$$\frac{\partial u_1}{\partial t} = \Delta u_1 + p(1 - u_1 - u_2) - r u_1 u_2,$$
$$\frac{\partial u_2}{\partial t} = \Delta u_2 + q(1 - u_1 - u_2)^2 - r u_1 u_2.$$

If $p < q$, the ODE has four fixed points: two stable $(1, 0)$ and $(\alpha, \beta)$ and two unstable: $(0, 1)$ and $(\beta, \alpha)$, where

$$\alpha < \beta = \frac{(q-p) \pm \sqrt{(q-p)^2 - 4qp^2/r}}{2q}.$$

See Figure 5 for an example.

The PDE results that were routine for the sexual reproduction model are now difficult. To prove the existence of a traveling wave with $u(-\infty) = (\alpha, \beta)$ and $u(\infty) = (1, 0)$ one goes to the four dimensional phase plane: $(u_1, u_1', u_2, u_2')$, and looks for a curve connecting $(\alpha, 0, \beta, 0)$ and $(1, 0, 0, 0)$ which will exist only for one value of the speed $c$. Fortunately this was done previously by Volpert and Volpert (1988). With the existence of a traveling wave established the next step is to prove a convergence theorem for the PDE, which can be done with comparison techniques because of monotonicity properties of system $(u_1, -u_2)$. Once the PDE result is established the rest is a routine application of the block construction.

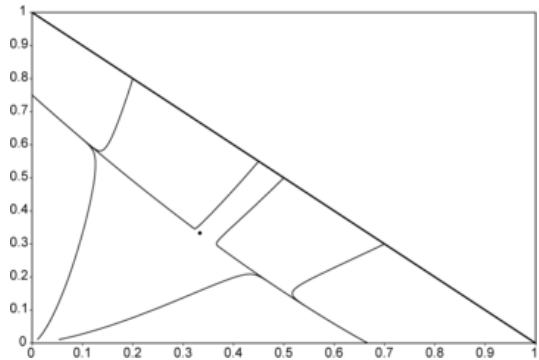

Fig. 6. *Colicin mean field ODE.*

*Example 2.3. Colicin.* Durrett and Levin (1997) considered a competition between two types of *E. coli*, one of which produces colicin (a chemical that kills other *E. coli*):

| Birth | Rate | Death | Rate |
|---|---|---|---|
| $0 \to 1$ | $\beta_1 f_1$ | $1 \to 0$ | $\delta_1$ |
| $0 \to 2$ | $\beta_2 f_2$ | $2 \to 0$ | $\delta_2 + \gamma f_1$ |

Here the rates are like the two-species contact process, except for the $\gamma f_1$ in the death rate $2 \to 0$, which comes from 1's killing 2's with colicin. For simplicity we suppose that the basic death rates are equal $\delta_1 = \delta_2 = 1$. Having done this it is natural to suppose that $\beta_1 < \beta_2$, for otherwise it is clear that the 1's will out compete the 2's.

The mean-field ODE is:

(9)
$$\frac{du_1}{dt} = \beta_1 u_1(1 - u_1 - u_2) - \delta_1 u_1,$$
$$\frac{du_2}{dt} = \beta_2 u_2(1 - u_1 - u_2) - u_2(\delta_2 + \gamma u_1).$$

A little algebra shows that when

$$\delta_i < \beta_i \quad \text{and} \quad \frac{\delta_2}{\beta_2} < \frac{\delta_1}{\beta_1} < \frac{\delta_2 + \gamma}{\beta_2 + \gamma}$$

the mean-field ODE has an interior fixed point but it is unstable. Figure 6 shows the situation when $\beta_1 = 3$, $\gamma_1 = 2.5$ and $\beta_2 = 4$.

Figure 7 gives the density versus time in the system on a $200 \times 200$ grid. To emphasize that the behavior is different from the ODE, we start the gray colicin producer (1's) at a small density. By time 1000 it has eliminated the black colicin sensitive strain (2's). The other panel shows the state at time 600. Note that the two types have segregated. A movie would show that the interface moves in a direction that favors the 1's.

PROBLEM 6. Show that coexistence is not possible in the colicin model when $\beta_1 < \beta_2$ and $\delta_1 = \delta_2 = 1$.

**Case 3. Cyclic systems, periodic orbits.** In this case, we see coexistence with significant spatial structure. The pictures are pretty but the problems are hard.

*Example 3.1. Multitype biased voter model.* Each site can be in state 1, 2, ..., $k$, and $j \to i$ at rate $f_i \lambda_{ij}$. In words, $i$'s eat $j$'s at rate $\lambda_{ij}$. The mean field ODE is

$$\frac{du_i}{dt} = u_i \sum_j (\lambda_{ij} - \lambda_{ji}) u_j.$$

Silvertown et al. (1992) who were interested in the competition of grass species, studied the five species case in which

$$\lambda_{ij} = \begin{pmatrix} 0 & 0.09 & 0.32 & 0.23 & 0.37 \\ 0.08 & 0 & 0.16 & 0.06 & 0.09 \\ 0.06 & 0.06 & 0 & 0.44 & 0.11 \\ 0.02 & 0.06 & 0.05 & 0 & 0.03 \\ 0.02 & 0.03 & 0.05 & 0.03 & 0 \end{pmatrix}.$$

This example is not very interesting because $\lambda_{1j} > \lambda_{j1}$ for $2 \leq j \leq 5$, so using ideas of Grannan and Swindle (1990) and their improvement by Mountford and Sudbury (1992) one can show that if $A_t^1$ is the event 1's are alive at time $t$ then $P(A_t^1, \xi_t(x) \neq 1) \to 0$ so if the 1's do not die out they take over the system. The key to the proof is that if $\theta$ is small

$$Z_t = \sum_{z:\xi_t(z)=1} e^{-\theta|z|} \qquad \text{is a submartingale.}$$

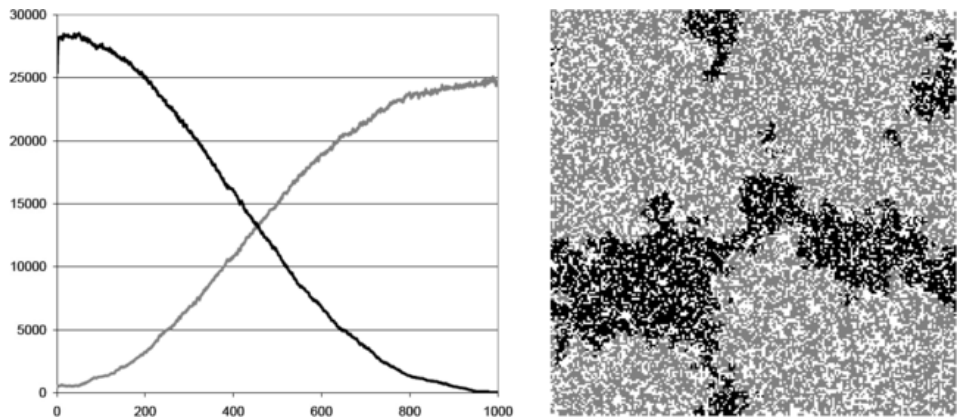

FIG. 7. *Simulation of colicin model. Producers (1's) are gray; sensitive strain (2's) is black.*

Durrett and Levin (1998) studied the cyclic case in which $\lambda_{13} = \beta_1$, $\lambda_{21} = \beta_2$, $\lambda_{32} = \beta_3$, and the other $\lambda_{ij} = 0$. This system with $\beta_i = 1$ and the corresponding discrete time deterministic cellular automata had been studied earlier by Bramson and Griffeath (1989), Fisch, Gravner and Griffeath (1991) and Durrett and Griffeath (1993). The mean-field ODE has equilibrium: $\rho_i = \beta_{i-1}/(\beta_1 + \beta_2 + \beta_3)$ where $i-1$ is computed modulo 3 with the result in $\{1,2,3\}$. Around this fixed point are concentric periodic orbits. See Figure 8 for an example. To prove mathematically that this occurs, write $H(u) = \sum_i \rho_i \log u_i$ and check that $H(u)$ is constant along solutions of the ODE.

PROBLEM 7. Show that coexistence always occurs in the cyclic case of Silvertown's model.

For partial credit show that coexistence can occur for some parameters. Figure 9 gives a proof by simulation.

*Interlude: lizard love.* Systems with a cyclic relationship exist in nature. In the side-blotched lizard (*Uta stansburiana*), males have one of three throat colors, each one declaring a particular strategy. Dominant, orange-throated males establish large territories within which live several females. But these territories are vulnerable to infiltration by males with yellow-striped throats—known as sneakers—who mimic the markings and behavior of receptive females. The orange males can't successfully defend all their females against these disguised interlopers, who cluster on the fringes of the territories held by the orange lizards.

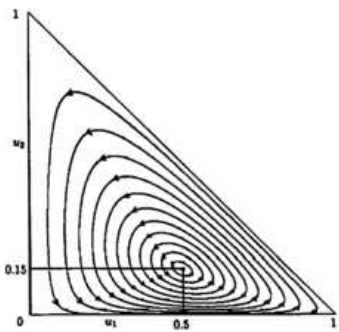

FIG. 8. *Cyclic particle system with $\beta_1 = 0.3$, $\beta_2 = 0.7$, $\beta_3 = 1.0$.*

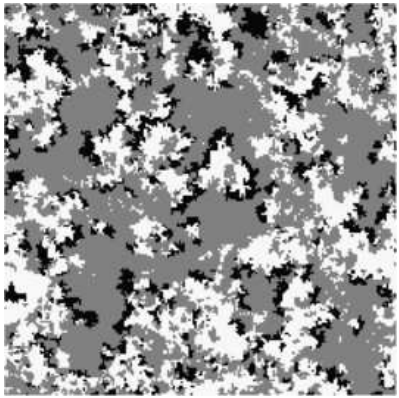

Fig. 9. *Cyclic particle system with $\beta_1 = 0.3$, $\beta_2 = 0.7$, $\beta_3 = 1.0$. 1's are black, 2's white, 3's gray.*

However, a large population of sneakers, which have no territory of their own to defend, can be quickly overrun by blue-throated males, who defend territories large enough to hold just one female. Sneakers have no chance against a vigilant, blue-throated guard. But once the sneakers become rare, powerful orange males flourish, grabbing territory and females from the blue lizards. Now, the blue males lose out. See Sinervo and Lively (1996) for more on these lizards.

*Example 3.2. Three species colicin.* Durrett and Levin (1997) considered an *E. coli* competition model with rates

| Birth | Rate | Death | Rate |
|---|---|---|---|
| $0 \to 1$ | $\beta_1 f_1$ | $1 \to 0$ | $\delta_1$ |
| $0 \to 2$ | $\beta_2 f_2$ | $2 \to 0$ | $\delta_2$ |
| $0 \to 3$ | $\beta_3 f_3$ | $3 \to 0$ | $\delta_3 + \gamma_1 f_1 + \gamma_2 f_2$ |

Here, 1's and 2's are colicin producers, while 3 is colicin sensitive. In the two species system (we conjecture) there is no coexistence, but as we will see coexistence is possible with three species.

Consider for concreteness, the situation when $\delta_i = 1$, $\beta_1 = 3$, $\beta_2 = 3.2$, $\beta_3 = 4.0$, $\gamma_1 = 3$ and $\gamma_2 = 0.5$. In this case the 2's beat the 1's since they have a larger birth rate, the 3's beat the 2's since the colicin they make is not nasty enough, while the 1's beat the 3's. Thus again the three competitors have the same relationship as in the child's game paper-rock-scissors.

The mean-field ODE is similar to (9):

$$\frac{du_1}{dt} = \beta_1 u_1 (1 - u_1 - u_2 - u_3) - \delta_1 u_1,$$

(10)
$$\frac{du_2}{dt} = \beta_2 u_2(1 - u_1 - u_2 - u_3) - \delta_2 u_2,$$
$$\frac{du_3}{dt} = \beta_3 u_3(1 - u_1 - u_2 - u_3) - u_3(\delta_3 + \gamma_1 u_1 + \gamma_2 u_2).$$

Figure 10 gives a picture of the mean-field ODE in the concrete case considered above as we look down into the tetrahedron $u_i \geq 0$, $u_1 + u_2 + u_3 \leq 1$. On the $(u_1, 0, u_3)$ and $(0, u_2, u_3)$ faces we see the colicin ODE, while on the $(u_1, u_2, 0)$ face we have the competing contact process. As in the competing contact process comparing the first two equations shows that there is fixed point with $u_1 u_2 > 0$ when $\beta_1/\delta_1 \neq \beta_2/\delta_2$.

Figure 11 gives a simulation. The graph gives the numbers of the different types on a $200 \times 200$ grid, while the picture gives a snapshot of part of the system at the final time.

PROBLEM 8. Show that coexistence can occur in the three species colicin model.

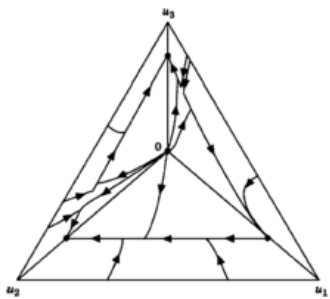

Fig. 10. *Three species colicin mean-field ODE.*

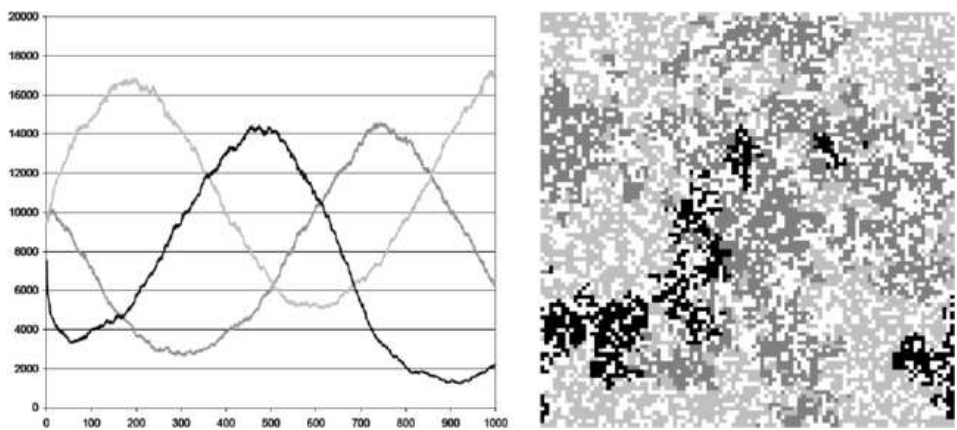

Fig. 11. *Three species colicin simulation: 1's black, 2's darker gray, 3's lighter gray.*

Coexistence has been verified experimentally by Kirkup and Riley (2004). They began with a sensitive strain ($S$) of *E. coli*, introduced colicin plasmids to make a colicin-producing strain ($C$), and exposed the sensitive strain to the colicin-producer to obtain a resistant strain ($R$). A number of unlucky mice were then chosen to have the competition drama play out in their colons. After reading the paper, I think I will stick to computer simulations. In four cases, the experiment had to be discontinued because the mice were fighting or several mice died.

*Example* 3.3. *Spatial Prisoner's Dilemma.* This time we allow multiple hawks $\eta_t(x)$ and doves $\zeta_t(x)$ at each site.

- *Migration.* Each individual at rate $\nu$ migrates to a nearest neighbor.
- *Death due to crowding.* Each individual at $x$ dies at rate $\kappa(\eta_t(x) + \zeta_t(x))$.
- *Game step.* Let $p_t(x)$ be the fraction of hawks in the $5 \times 5$ square centered at $x$. Hawks give birth (or death) at rate $ap_t(x) + b(1 - p_t(x))$, doves at rate $cp_t(x) + d(1 - p_t(x))$.

Here, "birth (or death)" means that if the quantity is positive it is a birth rate, but if it is negative it is $-1$ times a death rate.

An interesting choice for the game matrix is a Prisoner's Dilemma

|   | **H** | **D** |
|---|---|---|
| H | $a = -0.6$ | $b = 0.9$ |
| D | $c = -0.9$ | $d = 0.7$ |

The $H$ strategy dominates $D$, so it is the better choice, but the payoff for $(D, D)$ is better than that for $(H, H)$. This is the Prisoner's Dilemma "paradox." If everyone played $D$ then the world would be a nice place, but this leads to the temptation to play $H$ and increase your payoff.

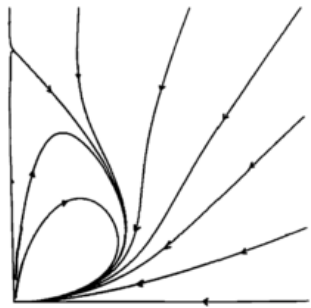

Fig. 12. *Hawks–Doves ODE.*

In a homogeneously mixing population the densities of Hawks ($u$) and Doves ($v$) would evolve according to

$$\frac{du}{dt} = u\left\{a\frac{u}{u+v} + b\frac{v}{u+v} - \kappa(u+v)\right\},$$
$$\frac{dv}{dt} = v\left\{c\frac{u}{u+v} + d\frac{v}{u+v} - \kappa(u+v)\right\}.$$
(11)

In the concrete example under consideration, we have the following behavior. On the vertical axis one can see that in the absence of Hawks, Doves reach an equilibrium. However, when both are present, the ratio of Hawks to Doves increases until the population crashes to 0. See Figure 12 for an example.

Simulations suggest that in the spatial model Hawks and Doves can coexist. Intuitively this occurs because the dynamics of the ODE happen locally, but when the Hawk population crashes to 0 then there are some Doves left behind to recolonize space, but when the Doves grow to a significant density then the remaining Hawks again take over.

PROBLEM 9. Show that coexistence can occur in the spatial Prisoner's Dilemma.

**Credits.** Since I first investigated the right edge of the one dimensional contact process in 1980, I have written 92 papers (out of my current total of 170) on the topic of stochastic spatial model. On 21 papers, I am the only author, but most of those are survey articles and conference proceedings, so much of the credit for my work should go to my collaborators. Listed in order of their multiplicity they are: (12) Ted Cox; (10) Simon Levin and Linda Buttel; (6) Maury Bramson; (4) David Griffeath, Roberto Schonmann; (3) Tom Liggett, Ed Perkins, Rinaldo Schinazi, Jeff Steif; (2) Wan-ding Ding, Nicolas Lanchier, Xiu-fang Liu, Iljana Zähle; (1) Ben Chan, Larry Gray, Paul Jung, Anne Moller, and Mateo Restrepo. In addition I would like to thank two referees, and my current postdoc John Mayberry for helping me to improve the paper.

## REFERENCES

BRAMSON, M. and DURRETT, R. (1988). A simple proof of the stability criterion of Gray and Griffeath. *Probab. Theory Related Fields* **80** 293–298. MR968822

BRAMSON, M. and GRIFFEATH, D. (1989). Flux and fixation in cyclic particle systems. *Ann. Probab.* **17** 26–45. MR972768

BRAMSON, M. and NEUAHUSER, C. (1992). A catalytic surface reaction model. *J. Comput. Appl. Math.* **40** 157–161.

CHAN, B. and DURRETT, R. (2006). A new coexistence result for competing contact processes. *Ann. Appl. Probab.* **16** 1155–1165. MR2260060

DURRETT, R. (1992). Multicolor particle systems with large threshold and range. *J. Theoret. Probab.* **5** 127–152. MR1144730

Durrett, R. (1995). Ten lectures on particle systems. In *Lectures on Probability Theory (Saint-Flour, 1993). Lecture Notes in Mathematics* **1608** 97–201. Springer, Berlin. MR1383122

Durrett, R. (2002). Mutual invadability implies coexistence in spatial models. *Mem. Amer. Math. Soc.* **156** viii–118. MR1879853

Durrett, R. and Griffeath, D. (1993). Asymptotic behavior of excitable cellular automata. *Experiment. Math.* **2** 183–208. MR1273408

Durrett, R. and Lanchier, N. (2008). Coexistence in host-pathogen systems. *Stochastic Process. Appl.* **118** 1004–1021. MR2418255

Durrett, R. and Levin, S. (1994). The importance of being discrete (and spatial). *Theor. Pop. Biol.* **46** 363–394.

Durrett, R. and Levin, S. (1997). Allelopathy in spatially distributed populations. *J. Theoret. Biol.* **185** 165–172.

Durrett, R. and Levin, S. (1998). Spatial aspects of interspecific competition. *Theor. Pop. Biol.* **53** 30–43.

Durrett, R. and Møller, A. M. (1991). Complete convergence theorem for a competition model. *Probab. Theory Related Fields* **88** 121–136. MR1094080

Durrett, R. and Neuhauser, C. (1994). Particle systems and reaction-diffusion equations. *Ann. Probab.* **22** 289–333. MR1258879

Durrett, R. and Neuhauser, C. (1997). Coexistence results for some competition models. *Ann. Appl. Probab.* **7** 10–45. MR1428748

Durrett, R. and Schinazi, R. (1993). Asymptotic critical value for a competition model. *Ann. Appl. Probab.* **3** 1047–1066. MR1241034

Durrett, R. and Swindle, G. (1991). Are there bushes in a forest? *Stochastic Process. Appl.* **37** 19–31. MR1091691

Durrett, R. and Swindle, G. (1994). Coexistence results for catalysts. *Probab. Theory Related Fields* **98** 489–515. MR1271107

Fisch, R., Gravner, J. and Griffeath, D. (1991). Cyclic cellular automata in two dimensions. In *Spatial Stochastic Processes. Progress in Probability* **19** 171–185. Birkhäuser Boston, Boston, MA. MR1144096

Grannan, E. R. and Swindle, G. (1990). Rigorous results on mathematical models of catalytic surfaces. *J. Statist. Phys.* **61** 1085–1103. MR1083897

Harris, T. E. (1974). Contact interactions on a lattice. *Ann. Probab.* **2** 969–988. MR0356292

Kirkup, B. C. and Riley, M. A. (2004). Antibiotic-meidated antagonism leads to a bacterial game of rock-paper-scissor in vivo. *Nature* **428** 412–414.

Lanchier, N. and Neuhauser, C. (2006). Stochastic spatial models of host-pathogen and host-mutualist interactions. I. *Ann. Appl. Probab.* **16** 448–474. MR2209349

Levin, S. A. (1970). Community equilibria and stability, and an extension of the competitive exclusion. principle. *Am. Naturalist.* **104** 413–423.

Liggett, T. M. (1985). *Interacting Particle Systems. Grundlehren der Mathematischen Wissenschaften [Fundamental Principles of Mathematical Sciences]* **276**. Springer, New York. MR776231

Liggett, T. M. (1999). *Stochastic Interacting Systems: Contact, Voter and Exclusion Processes. Grundlehren der Mathematischen Wissenschaften [Fundamental Principles of Mathematical Sciences]* **324**. Springer, Berlin. MR1717346

Mountford, T. S. and Sudbury, A. (1992). An extension of a result of Grannan and Swindle on the poisoning of catalytic surfaces. *J. Statist. Phys.* **67** 1219–1222. MR1170090

Neuhauser, C. (1992). Ergodic theorems for the multitype contact process. *Probab. Theory Related Fields* **91** 467–506. MR1151806

Silvertown, J., Holtier, S., Johnson, J. and Dale, P. (1992). Cellular automaton models of interspecific competition for space—the effect of pattern on process. *J. Ecol.* **80** 527–534.

Sinvervo, B. and Lively, C. M. (1996). The rock-paper-scissors game and the evolution of alternative male strategies. *Nature* **380** 240–243.

Volpert, V. A. and Volpert, A. I. (1988). Application of the Leray–Schauder method to the proof of the existence of wave solutions of parabolic systems. *Dokl. Akad. Nauk SSSR* **298** 784–787. MR936502

Ziff, R. M., Gulari, E. and Barshad, Y. (1986). Kinetic phase transitions in an irreversible surface-reaction model. *Phys. Rev. Lett.* **56** 2553–2556.

Department of Mathematics  
523 Malott Hall  
Cornell University  
Ithaca, New York 14853  
USA  
E-mail: rtd1@cornell.edu